\documentclass[reqno]{amsart}
\usepackage{hyperref}
\usepackage{amssymb}
\newtheorem{Lemma}{Lemma}[section]
\newtheorem{theorem}[Lemma]{Theorem}
\newtheorem{lemma}[Lemma]{Lemma}

\newtheorem{definition}[Lemma]{Definition}

\begin{document}
\title[On the uniqueness of the optimal control]
{On the uniqueness of the optimal control for two-dimensional
	second grade fluids}

\author[A. Almeida, N. V. Chemetov, F. Cipriano]
{Adilson Almeida, N. V. Chemetov, F. Cipriano}

\address{Adilson Almeida \newline
Centro de Matem\'atica e Aplica\c c\~oes (CMA), Faculdade de Ci\^encias e
Tecnologia da Universidade Nova de Lisboa, Portugal}
\email{ama.almeida@campus.fct.unl.pt}

\address{Nikolai V. Chemetov \newline
Department of Computing and Mathematics, University of S{\~a}o Paulo,
14040-901 Ribeir{\~a}o Preto - SP, Brazil}
\email{nvchemetov@gmail.com}

\address{Fernanda Cipriano \newline
Centro de Matem\'atica e Aplica\c c\~oes (CMA) FCT/UNL and Departamento de Matem\'atica, Faculdade de Ci\^encias e Tecnologia da Universidade Nova de Lisboa, Portugal}
\email{cipriano@fct.unl.pt}

\begin{abstract}
 We study an optimal control problem with a quadratic cost functional for  non-Newtonian fluids of differential type.
 More precisely, we consider the system governing the evolution of a second grade fluid filling a two-dimensional bounded domain, supplemented with a Navier slip boundary condition, and under certain assumptions on the size of the initial data and  parameters of the model,
 we  prove the second-order sufficient optimality conditions.
 Furthermore,  we establish a global uniqueness result for the solutions of the first-order optimality system. 
\end{abstract}

\maketitle

\section{Introduction}
The optimization of evolutionary phenomena is crucial in several branches of the knowledge, for instance in Finance, Biology, Ecology, Aviation etc.   
\cite{CMP21}, \cite{y1}, \cite{y2}, \cite{CC19}.   
The optimal control of fluid flows is a major problem in mathematical physics,
with relevant consequences in industrial applications. 
In the last decades, extensive research work has been carried out on the control of fluid flows described by the Navier-Stokes equations.
However, many incompressible viscous fluids present in the nature and used in the industry do not satisfy the Newton's law of viscosity, and consequently cannot be described by the Navier-Stokes equations.
Among these fluids,  called non-Newtonian fluids, we can  find
colloidal suspensions and emulsions, some industrial oils, ink-jet prints,
geological flows, biological fluids, body care fluids, some materials
arising in polymer processing as well as in food processing, and many others.

In this article, we study second grade fluids, which
belong to the class of non-Newtonian complex viscoelastic fluids of
differential type. To understand the physical principles associated to the second grade fluid equations, as well as the physical properties of these fluids, we refer to \cite{DF74}, \cite{DR95}, \cite{RE55} (see also
the more recent works \cite{AA11}, \cite{HAKA17}).

From the mathematical point of view, the equations governing the evolution
of second grade fluids are  strongly nonlinear partial differential
equations. The existence and uniqueness problems with a Dirichlet boundary
condition were established in the pioneering works \cite{CO84}, \cite{O81}
and \cite{CG97} (see also \cite{B20}). Despite the most usual boundary
condition to be the non-slip Dirichlet boundary condition, practical studies
show that some viscoelastic fluids slip against the boundary surface. Let us
refer, for instance, \cite{WD97} on capillary flow of highly entangled
polyethylene (PE) melts, and \cite{MBC04} on microgel pastes and
concentrated emulsions exhibiting a generic slip behavior at low stresses
when sheared near smooth surfaces. Therefore, to accurately describe certain
physical systems, a slip boundary condition should be considered (cf. \cite%
{CCG10}, \cite{CC_1_13}, \cite{CC_2_13}, \cite{CC_3_13}). The article \cite{BR03}
establishes the well-posedness for the second grade fluid equations under a
Navier slip boundary condition. Referring to the stochastic framework,
the existence and uniqueness results have been investigated in \cite{RS101}, 
\cite{RS12}, \cite{CC17}, \cite{S18}, \cite{SZZ19} under non-slip and slip
boundary conditions.

In the deterministic context, the optimal control problem for the second grade fluid equations was studied in the articles \cite{AC1}, \cite{AC2}, while
the optimal control of the stochastic dynamic has been addressed in \cite%
{CC18}, \cite{CP19}.  The authors established the existence of an optimal
solution for the control problem, and by analysing the linearized state
equation as well as the adjoint equation deduced the first-order optimality
conditions. 

In the present article, we perform a second-order analysis by studying the
second-order derivative of the objective functional, in addition, we obtain
a global uniqueness result for the optimal solution.
 Essentially, if the fluid material is filling a bounded two-dimensional domain, the sufficient second-order optimality condition is achieved for sufficiently elastic and viscous fluid material with small size initial
data. 
Alternatively, the same result can be obtained by taking an objective
functional with strong intensity of the cost.

We should emphasize that the second-order sufficient optimality conditions
are crucial to numerically solve the control problem, being necessary for
the stability of the optimal solution, as well as to prove the
convergence of the numerical approximations.
 As far as we know, the second-order analysis and the uniqueness problem for the optimization of second
grade fluids is being addressed here for the first time.

In Section 2, we formulate the problem and recall some known results in the
literature. We write key estimates for the solutions of the state equation,
linearized equation and adjoint equation that will be necessary in the
following sections. Section 3 establishes the  second-order sufficient optimality conditions. Section 5 is devoted to the global uniqueness problem.

\section{Formulation of the problem and preliminary results}

\bigskip \setcounter{equation}{0}\label{5}

We consider an optimal control problem associated with a non-stationary
viscous, incompressible, second grade fluid. We assume that $\mathcal{O}%
\subset \mathbb{R}^{2}$ is a bounded, simply connected domain, having a
sufficiently smooth boundary $\Gamma .$ The fluid dynamic on a time
interval $[0,T]$ is described by the following state equations 
\begin{equation}
	\left\{ 
	\begin{array}{cc}
		\frac{\partial }{\partial t}\upsilon (y)=\nu \Delta y-\mathrm{curl}%
		\,\upsilon (y)\times y-\nabla \pi +u, & \vspace{1mm} \\ 
		\multicolumn{1}{l}{\mathrm{div}\,y=0} & \multicolumn{1}{l}{\mbox{in}\
			Q=(0,T)\times \mathcal{O},} \\ 
		\multicolumn{1}{l}{y\cdot \mathrm{n}=0,\qquad (\mathrm{n}\cdot Dy)\cdot 
			\mathrm{\tau }=0} & \multicolumn{1}{l}{\mbox{on}\ \Sigma =(0,T)\times \Gamma
			,} \\ 
		\multicolumn{1}{l}{y(0)=y_{0}} & \multicolumn{1}{l}{\mbox{in}\ \mathcal{O},}%
	\end{array}%
	\right.  \label{equation_etat}
\end{equation}%
where 
$y=(y_{1},y_{2})$ is the velocity field of the fluid, $$Dy=\frac{ \nabla y+(\nabla
y)^{T} }{2}$$ corresponds to the symmetric part of the velocity gradient and 
\begin{equation*}
	\upsilon (y)=y-\alpha \Delta y,
\end{equation*}%
where $\alpha >0$ is a viscoelastic parameter. 
Moreover $\pi $ denotes the
hydrodynamic pressure, \ $\nu >0$ is the viscosity of the fluid, $\mathrm{n}%
=(n_{1},n_{2})$ is the unit normal to the boundary $\Gamma $, $\mathrm{\tau }%
=(-n_{2},n_{1})$ is the tangent vector to $\Gamma $ and $u$ represents an
external mechanical force, which acts on the system as the control variable.
Here $\upsilon=\upsilon (y)$ and $y$ are two-dimensional vectors.  To perform a three-dimensional calculus in the term $\mathrm{curl}%
\,\upsilon (y)\times y$ (or in a similar ones) we consider the usual identifications  $$\mathrm{curl}%
\,v=\tfrac{\partial v_{2}}{\partial x_{1}}-\tfrac{\partial v_{1}}{\partial
	x_{2}}\equiv(0,0, \tfrac{\partial v_{2}}{\partial x_{1}}-\tfrac{\partial v_{1}}{\partial
	x_{2}})\quad \text{and}\quad  y=(y_1,y_2)\equiv (y_1,y_2,0).$$ 

\bigskip

In order to formulate the problem and establish the results, we introduce
the convenient functional spaces and some useful notations.

We denote by $L^{p}=L^{p}(\mathcal{O})$, $1 \leq p \leq \infty $, the  Lebesgue spaces, endowed  with their natural norms $%
\Vert \cdot \Vert _{p}$ and  consider the usual notations for the scalar products on the finite
dimensional spaces $\mathbb{R}^{2}$ and $\mathbb{R}^{2\times 2}$ 
\begin{equation*}
	u\cdot v=\sum_{i=1}^{2}u_{i}v_{i},\quad u,v\in \mathbb{R}^{2};\quad \eta
	:\zeta =\sum_{i,j=1}^{2}\eta _{ij}\zeta _{ij},\quad \eta ,\zeta \in \mathbb{R%
	}^{2\times 2},
\end{equation*}%
as well as for the scalar products in $L^{2}$ 
\begin{equation*}
	\left( u,v\right) =\displaystyle\int_{\mathcal{O}}u(x)\cdot v(x)\,dx,\qquad
	\left( \eta ,\zeta \right) =\displaystyle\int_{\mathcal{O}}\eta (x):\zeta
	(x)\,dx.
\end{equation*}%
We consider the standard  Sobolev spaces $W^{k,p}=W^{k,p}(\mathcal{O})$, $1 \leq p \leq \infty $,  endowed  with their natural norms $\Vert \cdot \Vert _{W^{k,p}}$. In the particular case $p=2$, we set 
$H^k=W^{k,2}$ and 
 $\Vert
\cdot \Vert _{H^{k}}=\Vert \cdot \Vert _{W^{k,2}}$.

\bigskip

Let us introduce the
following divergence-free spaces 
\begin{equation*}
	\begin{array}{ll}
		H & =\left\{ v\in L^{2}\mid \mathrm{div}\,v=0\ \text{ in } \mathcal{D}^\prime
		(\mathcal{O})\ \mbox{ and }\ v\cdot n=0\ \mbox{ on }H^{-1/2}(\Gamma
		)\right\} ,\vspace{2mm} \\ 
		V & =\left\{ v\in H^{1}\mid \mathrm{div}\,v=0\ \mbox{ in }\ 
		\mathcal{O}\mbox{ and }\ v\cdot n=0\ \text{ on }\ \Gamma \right\} ,\vspace{%
			2mm} \\ 
		W & =\left\{ v\in V\cap H^{2}\mid \left( n\cdot Dv\right) \cdot
		\tau =0\ \ \mbox{on}\ \Gamma \right\} .%
	\end{array}%
\end{equation*}%
On the space $V,$ we consider the following scalar product and the corresponding norm%
\begin{equation*}
	\left( u,v\right) _{V}=\left( u,v\right) +2\alpha \left( Du,Dv\right)
	,\qquad \Vert u\Vert _{V}=\left( u,u\right) _{V}^{1/2},\qquad \forall u,v\in
	V.
\end{equation*}

Throughout the article,  we denote 
by $C,C_{1},...,C_{4}$ \ and $K$, $\widetilde{K},$ $\widehat{K}$,
the constants which depend only on the domain $\mathcal{O}.$

 Here, we
recall the usual Korn inequality 
\begin{equation}
	\Vert y\Vert _{H^{1}}^{2}\leq K\left( \Vert y\Vert _{2}^{2}+\Vert Dy\Vert
	_{2}^{2}\right) ,\quad \forall y\in H^{1}.  \label{korn_1}
\end{equation}%
Let us define the operator 
$$\mathbb{A}y\mathbb{=}\mathbb{P\triangle }y\quad 
\text{for} \quad y\in W,$$
 where $\ \mathbb{P}:L^{2}\rightarrow H$ is the
Helmholtz projector in $L^{2}.$ 
Due to 
 Lemma 2.3 in \cite{AC1},  we have the following inequality 
\begin{equation}
	\Vert y\Vert _{H^{2}}^{2}\leq \widetilde{K}\left( \Vert y\Vert
	_{2}^{2}+\Vert \mathbb{A}y\Vert _{2}^{2}\right) ,\quad \forall y\in W\cap
	H^{3}.  \label{a}
\end{equation}%

\bigskip

Considering the trilinear form 
\begin{equation*}
	b\left( \phi ,z,y)=(\phi \cdot \nabla z,y\right) ,
\end{equation*}%
the nonlinear term of the equations can be written as 
\begin{equation}
	\begin{array}{ll}
	\left(curl\,\upsilon (y)\times z,\phi \right) =b\left( \phi
,z,\upsilon (y)\right) -b\left( z,\phi ,\upsilon (y)\right), \quad  \forall
y,z\in W\cap H^{3},\;  \phi \in V. \label{as} 
	\end{array}%
\end{equation}%

Now, let us state a property of the nonlinear term, which will be useful in the following considerations.
\begin{lemma}
	\label{L1} For all $z,$ $\phi $ $\in H^{2}$, we have  
	\begin{equation}
		\left\vert \left( \mathrm{curl}\,\upsilon (z)\times z,\phi \right)
		\right\vert \leq \widehat{K}\Vert \phi \Vert _{H^{2}}\Vert z\Vert
		_{H^{2}}^{2}.  \label{21:49:25S}
	\end{equation}
\end{lemma}

\textbf{Proof. } We verify that
\begin{align*}
	\left\vert \left( \mathrm{curl}\,\upsilon (z)\times z,\phi \right)
	\right\vert & \leq |b(\phi ,z,\upsilon (z))|+|b(z,\phi ,\upsilon (z))| \\
	& \leq \Vert \phi \Vert _{4}\Vert z\Vert _{W^{1,4}}\Vert z\Vert
	_{H^{2}}+\Vert z\Vert _{4}\Vert \phi \Vert _{W^{1,4}}\Vert z\Vert
	_{H^{2}}\leq \widehat{K}\Vert \phi \Vert _{H^{2}}\Vert z\Vert _{H^{2}}^{2}.
\end{align*}
\hfill $\blacksquare $
\bigskip

In this article, we assume  that the external mechanical force $u$ and the initial data $y_{0}$ satisfy

\begin{equation}
	u\in L^{2}(0,T;H^{1})\quad \text{and}\quad y_{0}\in W\cap H^{3}.
	\label{data}
\end{equation}

\begin{lemma}
	\label{L2} Under the assumptions \eqref{data} there exists a unique solution 
	$$
	y\in L^{\infty }\left( 0,T;W\cap H^{3}\right) \text{ with } \frac{%
		\partial y(t)}{\partial t}\in L^{2}\left( 0,T;V\right) 
	$$ 
	of the problem $(\ref{equation_etat}),$ which is understood in the distributional sense
	\begin{align}
		\big(\frac{\partial y(t)}{\partial t},\phi \bigr)& +2\alpha \bigl(D\frac{%
			\partial y(t)}{\partial t},D\phi \bigr)+2\nu \left( Dy(t),D\phi \right)  
		\notag \\
		& +\left( \mathrm{curl}\,\upsilon (y(t))\times y(t),\phi \right) =\left(
		u(t),\phi \right) ,\;\;\forall \phi \in V.  \label{var_form_state}
	\end{align}%
	Moreover, the solution $y$ verifies the following  estimate 
	\begin{equation}
		\left\Vert y\right\Vert _{L^{\infty }{(0,T;H^{3})}}^{2}   \leq \frac{%
			C_{1}^{2}\lambda _{1}^{2}}{\alpha ^{2}} \label{12:12:20:10:10}
	\end{equation}%
	with $$ \lambda _{1}^{2}=\left( 1+4K\widehat{\alpha }\right)
	(\left\Vert y_{0}\right\Vert _{H^{3}}^{2}+\left\Vert u\right\Vert _{L^{1}{%
			(0,T;H^{1})}}^{2}),$$ where the constant $K$ is defined by \eqref{korn_1},
	and $\widehat{\alpha }=\max (\left( 2\alpha \right) ^{-1},2\alpha
	).$
\end{lemma}

\bigskip \textbf{Proof. }The solvability of the problem $(\ref{equation_etat}%
)$ is shown in  \cite{BR03}. Let us write in a convenient form the  estimates for the state variable $y$ with
respect to the initial data $y_{0}$ and the control variable $u$.  The estimates (3.4) and (3.5) 
of Proposition 3.2 in \cite{AC1}\ give  
\begin{equation}
\label{03:03:1}
	\Vert y\Vert _{L^{\infty }{(0,T;H^{3})}}^{2}\leq \frac{C^{2}}{\alpha
		^{2}}\left( \Vert y\Vert _{L^{\infty }{(0,T;H^{1})}}^{2}+\Vert \mathrm{curl}%
	\sigma (y_{0})\Vert _{{L^{2}}}^{2}+\Vert \mathrm{curl}u\Vert _{L^{1}{%
			(0,T;L^{2})}}^{2}\right) .
\end{equation}%
In addition, the relation (3.3) of Proposition 3.2 in \cite{AC1} for $y$ in $H^{1},$ and the Korn's
inequality \eqref{korn_1} yield%
\begin{equation}
\label{03:03:2}
	\left\Vert y\right\Vert _{L^{\infty }{(0,T;H^{1})}}^{2}\leq 4K\frac{\max
		\left( 1,2\alpha \right) }{\min \left( 1,2\alpha \right) }\left( \left\Vert
	y_{0}\right\Vert _{H^{1}}^{2}+\left\Vert u\right\Vert
	_{L^{1}(0,T;L^{2})}^{2}\right) .
\end{equation}%
Combining the inequalities \eqref{03:03:1}-\eqref{03:03:2} we derive \eqref{12:12:20:10:10}.\hfill $%
\blacksquare $



\bigskip
Let us introduce the so-called \textit{control-to-state} mapping $%
S:u\rightarrow y,$ namely $y=S(u)$ is  the solution of the equation $ (\ref{equation_etat}%
)$ corresponding to the control $u$, and consider the reference pair $(u,y=S(u))$. In the next lemma, we recall a
stability result for the solution of (\ref{equation_etat}), which was
proved in \cite{AC1}, Proposition 4.4.

\begin{lemma}
	\label{L3}Let us consider the initial data $y_{0}$ and two different control variables  $%
	u_{1},u_{2}$\ satisfying the assumptions \eqref{data}. Let $y_{1}=S(u_{1})$, $%
	y_{2}=S(u_{2})$ be the corresponding  solutions of $(\ref{equation_etat})$ with
	the same initial data $y_{0}$. Then the difference $\overline{y}=y_{2}-y_{1}$
	verifies the estimate 
	\begin{equation}
		\Vert \overline{y}\Vert _{L^{\infty }(0,T;H^{2})}^{2}\leq \lambda
		_{2}^{2}\Vert \overline{u}\Vert _{L^{2}{(Q)}}^{2}  \label{est2}
	\end{equation}%
	with  $\overline{u}=u_{2}-u_{1}$ and 
	\begin{equation*}
		\lambda _{2}^{2}=\widetilde{K}\left[ \left( 1+C_{2}T(1+\alpha ^{-1})\frac{%
			C_{1}\lambda _{1}}{\alpha }\right) e^{C_{2}T\left( 1+\left( 1+\alpha
			^{-1}\right) C_{1}\lambda _{1}\right) }+\left( \alpha \nu \right) ^{-1}%
		\right] 
	\end{equation*}%
	for some positive constant $C_{1},C_{2}$ depending only on $\mathcal{O}.$
\end{lemma}

\textbf{Proof. }By the first estimate of Proposition 4.4 in \cite{AC1} and %
\eqref{12:12:20:10:10},\ we deduce 
\begin{equation*}
	\Vert \overline{y}\Vert _{L^{\infty }{(0,T;L^{2})}}^{2}+2\alpha \Vert D%
	\overline{y}\Vert _{L^{\infty }{(0,T;H^{1})}}^{2}\leq  e^{C_{2}T\left(
		1+\left( 1+\alpha ^{-1}\right) C_{1}\lambda _{1}\right) }\Vert \overline{u}%
	\Vert _{L^{2}{(Q)}}^{2}.
\end{equation*}%
Using the second estimate of Proposition 4.4 in \cite{AC1} and \eqref{a}, we
derive 
\begin{eqnarray*}
	\Vert \overline{y}\Vert _{L^{\infty }{(0,T;H^{2})}}^{2} &\leq &	\widetilde{K} \left[ \Vert \overline{y}\Vert _{L^{\infty }{(0,T;L^{2})}}^{2}+\right. \\
	&& \left. \frac{1}{\alpha }\left( \frac{1}{\nu }\Vert \overline{u}\Vert _{L^{2}{(Q)}%
	}^{2}+C_{2}T(1+\alpha ^{-1})C_{1}\lambda _{1}\Vert \overline{y}\Vert
	_{L^{\infty }{(0,T;L^{2})}}^{2}\right) \right]  \\
	&\leq &\widetilde{K} \left[ \left( 1+C_{2}T(1+\alpha ^{-1})\frac{%
		C_{1}\lambda _{1}}{\alpha }\right) \Vert \overline{y}\Vert _{L^{\infty }{%
			(0,T;L^{2})}}^{2}+ \right.\\
		&& \left. \left( \alpha \nu \right) ^{-1}\Vert \overline{u}\Vert
	_{L^{2}{(Q)}}^{2} \right] . \notag
\end{eqnarray*}%
These two inequalities imply \eqref{est2}.\hfill $\blacksquare $

\bigskip

According to  Proposition 4.5 in \cite{AC1}, the first-order G\^{a}teaux
derivative $$z=S^{\prime }(u)[w]$$
 of the mapping $S$, 
at the point $u$, in the
direction $w,$ is given by the  solution of the following linearized state equation
at $(u,y)$ 
\begin{equation}
	\left\{ 
	\begin{array}{ll}
		\displaystyle\frac{\partial \upsilon (z)}{\partial t}-\nu \Delta z+\mathrm{%
			curl}\,\upsilon (z)\times y+\mathrm{curl}\,\upsilon (y)\times z+\nabla \pi
		=w, &  \\ 
		\nabla \cdot z=0 & \mbox{in}\ Q, \\ 
		z\cdot n=0,\quad (n\cdot Dz)\cdot \tau =0 & \mbox{on}\ \Sigma , \\ 
		z(0)=0 & \mbox{in}\ \mathcal{O},%
	\end{array}%
	\right.  \label{linearized}
\end{equation}%
which   is well-posed in the Sobolev space $%
H^{2}$.

The next lemma establishes a suitable estimate for the solution $z.$ 
\begin{lemma}
	\label{L4} Under the assumptions \eqref{data} there exists a unique
	solution $z\in L^{\infty }(0,T;W)$ of \eqref{linearized}, such that 
	\begin{equation}
		\Vert z\Vert _{L^{\infty }(0,T;H^{2})}^{2}\leq \lambda _{3}^{2}\left\Vert
		w\right\Vert _{L^{2}(Q)}^{2},  \label{12:12:20:09:57}
	\end{equation}
	\\
	where 
\begin{align*}
\lambda _{3}^{2} =&\widetilde{K}\biggl[ \left( \alpha \nu \right) ^{-1}e^{C_{3}TC_{1}\lambda _{1}\left( \frac{1+\alpha }{\alpha ^{2}}\right)
}\biggl( 1+\frac{2}{\alpha }C_{3}TC_{1}\lambda _{1}\left( 1+\alpha
^{-1}\right) \text{$e$}^{C_{3}TC_{1}\lambda _{1}\left( \frac{1+\alpha }{%
\alpha ^{2}}\right) }\biggr) \times \\
&\qquad   e^{C_{3}T\left( 1+C_{1}\lambda _{1}
\left(1+\alpha ^{-1}\right) \right) }\biggr]
\end{align*}
	with $\lambda _{1}$  defined in \eqref{12:12:20:10:10}.
\end{lemma}

\textbf{Proof.} From the first estimate of Proposition 4.3 in \cite{AC1} and %
\eqref{12:12:20:10:10}, we have 
\begin{equation}
	\Vert z\Vert _{L^{\infty }(0,T;L^{2})}^{2}\leq \Vert w\Vert
	_{L^{2}(Q)}^{2}e^{C_{3}T\left( 1+C_{1}\lambda _{1}\left( 1+\alpha
		^{-1}\right) \right) }.  \label{aaa}
\end{equation}%
In addition, the second estimate of Proposition 4.3 in \cite{AC1} and %
\eqref{12:12:20:10:10} yield
\begin{eqnarray} 
\Vert \mathbb{A}z\Vert _{L^{\infty }(0,T;L^{2})}^{2}&\leq& \frac{1}{\alpha }%
\left[ \frac{1}{\nu }\left\Vert w\right\Vert _{L^{2}{(Q)}%
}^{2}+2C_{3}TC_{1}\lambda _{1}\left( 1+\alpha ^{-1}\right) \left\Vert
z\right\Vert _{L^{\infty }(0,T;L^{2})}^{2}\right] \times \notag
\\
&& e
^{C_{3}TC_{1}\lambda _{1}\left( \frac{1+\alpha }{\alpha ^{2}}\right) }=R. \label{03:03:4}
\end{eqnarray}%
Therefore, taking into account \eqref{a} and \eqref{03:03:4}, we deduce
\begin{eqnarray*}
	\Vert z\Vert _{L^{\infty }{(0,T;H^{2})}}^{2} &\leq &\widetilde{K}\left[
	\Vert z\Vert _{L^{\infty }{(0,T;L^{2})}}^{2}+R\right]  \\
	&=&\widetilde{K}\left[ \left( 1+\frac{1}{\alpha }2C_{3}TC_{1}\lambda
	_{1}\left( 1+\alpha ^{-1}\right) \text{$e$}^{C_{3}TC_{1}\lambda _{1}\left( 
		\frac{1+\alpha }{\alpha ^{2}}\right) }\right) \times \right.  \\
	&&\left. \Vert z\Vert _{L^{\infty }{%
			(0,T;L^{2})}}^{2}+\left( \alpha \nu \right) ^{-1}\text{$e$}^{C_{3}TC_{1}\lambda
		_{1}\left( \frac{1+\alpha }{\alpha ^{2}}\right) }\Vert w\Vert _{L^{2}{(Q)}%
	}^{2}\right] .
\end{eqnarray*}%
Using \eqref{aaa}, we obtain \eqref{12:12:20:09:57}.\hfill $\blacksquare $


\bigskip

Now, we formulate the control problem (see  \cite{AC1}). The space $\mathcal{U}_{ad}$ of the admissible control variables   is a bounded, closed and convex subset of $%
L^{2}(0,T;H^{1}).$ Namely, there exists a constant 
$L>0$ such that 
\begin{equation}
	\Vert u\Vert _{L^{2}(0,T;H^{1})}\leq L,\qquad \forall u\in 
	\mathcal{U}_{ad}.  \label{Uad}
\end{equation}%
The control on the evolution of the physical system is imposed through a distributed mechanical force $u\in 
\mathcal{U}_{ad}$, aiming to match a desired target velocity profile 
\begin{equation*}
	y_{d}\in L^{2}(Q).
\end{equation*}%
The control $u$ and the state $y=S(u)$ are constrained to
satisfy the system \eqref{equation_etat}, and the optimal control problem reads
\begin{equation*}
	(\mathcal{P})\qquad \qquad \underset{u}{\mbox{minimize}}\left\{ J(u,y):~u\in 
	\mathcal{U}_{ad},\quad y=S(u)\right\}, 
\end{equation*}%
where the objective functional is defined by
\begin{equation}
	\displaystyle J(u,y)=\;\frac{1}{2}\int_{0}^{T}\int_{\mathcal{O}}\left\vert
	y-y_{d}\right\vert ^{2}\,dxdt+\frac{\lambda }{2}\int_{0}^{T}\int_{\mathcal{O}%
	}\left\vert u\right\vert ^{2}\,dxdt,\qquad u\in \mathcal{U}_{ad}^{b},
	\label{cost}
\end{equation}%
and  $\lambda \geq 0$ is a fixed cost coefficient.

\bigskip

Let us recall that  the first-order
optimality conditions, at a locally optimal pair $(u^{\ast },y^{\ast })$, can be formally deduced through  the Lagrange's multipliers
method. More precisely, considering the Lagrange function 
\begin{equation*}
	L(u,y,p)=J(u,y)-\int_{0}^{T}\left( p,-\frac{\partial }{\partial t}\upsilon
	(y)+\nu \Delta y-\mathrm{curl}\,\upsilon (y)\times y-\nabla \pi +u\,\right)
	dt,
\end{equation*}%
the optimal pair $(u^*,y^*)$ should satisfy the following relations 
\begin{align}
	& L_{y}^{\prime }(u^{\ast },y^{\ast },p)[z]=0,\vspace{2mm}
	\label{09:46:27:10:20} \\
	& L_{u}^{\prime }(u^{\ast },y^{\ast },p)[w]\geq 0,\quad \forall w\in 
	\mathcal{U}_{ad}.  \label{20:42:14:11:20}
\end{align}%
The equation \eqref{09:46:27:10:20} yields the adjoint equation of the state
equation linearized at the point $y=y^{\ast }$, and the system composed by
the equations \eqref{09:46:27:10:20}, \eqref{20:42:14:11:20} and the state
equations \eqref{equation_etat} for $(u^{\ast },y^{\ast })$ constitute the so-called first-order
optimality conditions.

\bigskip

The adjoint equation of the state equation linearized at the point $%
y=y^{\ast }$ arising from \eqref{09:46:27:10:20} is given by 
\begin{equation}
	\left\{ 
	\begin{array}{ll}
		\displaystyle\frac{\partial }{\partial t}\upsilon (p)=-\nu \Delta p-\mathrm{%
			curl}\,\upsilon (y)\times p+\mathrm{curl}\,\upsilon (y\times p)+\nabla {\pi }%
		-(y-y_{d}),  \\ 
		\mathrm{div}\,p=0 \qquad\qquad\qquad\qquad\qquad\qquad\qquad \qquad\qquad\qquad \mbox{in}\ Q, \\ 
		p\cdot \mathrm{n}=0,\qquad\qquad (\mathrm{n}\cdot Dp)\cdot \mathrm{\tau }=0 \qquad \qquad\quad\qquad\qquad \,
		\mbox{on}\ \Sigma , \\ 
		p(T)=0\qquad\qquad\qquad\qquad\qquad\qquad\qquad \qquad\qquad\qquad \mbox{in}\ \mathcal{O}.%
	\end{array}%
	\right.  \label{adj_opt_eq_alpha}
\end{equation}%

Here, we recall the main result in \cite{AC1}, which shows the existence of an optimal solution and establishes  the  first-order optimality conditions.
\begin{theorem}
	\label{main_existence}\cite{AC1} Under the assumptions \eqref{data} and \eqref{Uad}, 
	the problem $(\mathcal{P})$ admits at least one solution 
	\begin{equation*}
		(u,y)=(u,S(u))\in \mathcal{U}_{ad}\times  L^{\infty }(0,T;W\cap H^{3}).
	\end{equation*}%
	Furthermore, there exists a unique solution 
	\begin{equation*}
		p=S^{\ast }(u)\in L^{\infty }( 0,T;H^{2} ) 
	\end{equation*}%
	of the adjoint system \eqref{adj_opt_eq_alpha}, such that the following
	optimality condition holds 
	\begin{equation}
		\int_{0}^{T}(\lambda u+p,\Psi -u)\,dt\geq 0,\qquad \forall \Psi \in \mathcal{%
			U}_{ad}.  \label{opt1}
	\end{equation}
\end{theorem}

A triplet  $(u^*, y^*, p^*)$ obtained as a solution   of the  coupled system (constituted by 
the state equation \eqref{equation_etat}, the adjoint equation \eqref{adj_opt_eq_alpha} and  the variational inequality \eqref{opt1}) 
is a candidate for an optimal solution, but not necessarily an optimal solution.
   The goal of this article is to 
analyse the conditions (on the initial data or on the parameters of the model) which  guarantee that the solution  of the  coupled system is optimal and unique.
This is a crucial step towards the implementation of the numerical methods to approximate the optimal control.

To perform this task,  we start by establishing convenient estimates for the adjoint state $p$. 

\begin{lemma}
	\label{L5} Under the assumptions of Theorem \ref{main_existence}, the adjoint
	state $p$ fulfils the following property
	\begin{align}
	&&	\Vert p\Vert _{L^{\infty }(0,T;H^{2})}^{2}\leq  \lambda _{4}^{2}=2%
		\widetilde{K}\left[ \left( 1+\frac{1}{\alpha }C_{4}TC_{1}\lambda _{1}\left(
		1+\alpha ^{-1}\right) \text{$e$}^{C_{4}TC_{1}\lambda _{1}\left( \frac{%
				1+\alpha }{\alpha ^{2}}\right) }\right) \times \right.
			\notag \\
		&& \left. e^{C_{4}T\left( 1+C_{1}\lambda
			_{1}\left( 1+\alpha ^{-1}\right) \right) }+
		 \left( \alpha \nu \right) ^{-1}\text{$e$}^{C_{4}TC_{1}\lambda _{1}%
			\frac{\left( 1+\alpha \right) ^{2}}{\alpha }}\right] \left( \frac{%
			C_{1}^{2}\lambda _{1}^{2}}{\alpha ^{2}}+\left\Vert y_{d}\right\Vert _{L^{2}{%
				(Q)}}^{2}\right) ,  \label{13:12:20:11:24}
	\end{align}%
where the constants $\widetilde{K}$ and $\lambda _{1}$
	are defined by \eqref{a} and \eqref{12:12:20:10:10}, respectively.
\end{lemma}

\textbf{Proof. } Taking into account  the first estimate of Proposition 5.4
in \cite{AC1}, we have 
\begin{equation}
\label{03:03:5}
	\Vert p\Vert _{L^{\infty }(0,T;L^{2})}^{2}\leq \Vert y-y_{d}\Vert _{L^{2}{(Q)%
	}}^{2}e^{C_{4}T\left( 1+C_{1}\lambda _{1}\left( 1+\alpha ^{-1}\right)
		\right) }.
\end{equation}
The second estimate of the same proposition gives 
	\begin{align}
	\label{03:03:6}
		\Vert \mathbb{A}p\Vert _{L^{\infty }(0,T;L^{2})}^{2}&\leq
 \frac{1}{\alpha }%
	\left( \tfrac{1}{\nu }\left\Vert y-y_{d}\right\Vert _{L^{2}{(Q)}%
	}^{2}+C_{4}TC_{1}\lambda _{1}\left( 1+\alpha ^{-1}\right) \left\Vert
	p\right\Vert _{L^{\infty }(0,T;L^{2})}^{2}\right) \times	\notag \\
	& \qquad e^{C_{4}TC_{1}\lambda _{1}\frac{\left( 1+\alpha \right) ^{2}}{\alpha }}.
\end{align}

Therefore, using the inequalities
 \eqref{03:03:5}-\eqref{03:03:6} and \eqref{a}, we can  apply 
the same reasoning as in the proof of Lemma \ref{L4} to  obtain   the claimed result.
\hfill $\blacksquare $

\bigskip

\section{Uniqueness results for the control problem $(\mathcal{P})$}

\bigskip \setcounter{equation}{0}\label{6}

This section establishes the main results of the article. First, by analysing the second-order derivative of the control-to-state mapping, as well as the second-order derivative of the cost functional, we deduce  a sufficient second-order
optimality condition, which guarantees that any triplet $(u^*,y^*, w)$ obtained as a solution of the first-order coupled optimality system will produce a locally optimal pair  $(u^*,y^*)$.

Next, we will be able to prove that the solution of the first-order coupled optimality system is unique. This uniqueness result 
conjugated with the result of Theorem \ref{main_existence} yields
the uniqueness of the solution for the optimal control problem $(\mathcal{P})$.
Our results will be achieved 
under some natural assumptions relying on the size of the initial data and  the parameters of the model, or  on the intensity of the cost. 
Essentially, if the fluid material is sufficiently viscous and elastic and the initial condition is small enough, or instead if the intensity of the cost is big enough,
the solution of the first-order optimality system is unique, and corresponds to the unique global solution of the optimal control problem $(\mathcal{P})$.

\bigskip

\subsection{Sufficient second-order optimality conditions}

In this section, we perform a second-order analysis. Let $y=S(u)$ be the
solution of $(\ref{equation_etat})$ for $u$ (and $y_{0}$)  verifying %
\eqref{data}. Let us consider some $w_{i}$, $i=1,2,$ satisfying \eqref{data}
too. Denoting $$z_{i}=S^{\prime }(u)[w_{i}], \quad i=1,2,$$ we can use standard
arguments to show that the second-order G\^{a}teaux derivative 
\begin{equation*}
	\widetilde{z}=S^{\prime \prime }(u)[w_{1},w_{2}]=\frac{d}{d\epsilon }\bigg|%
	_{\epsilon =0}S^{\prime }(u+\epsilon w_{2})[w_{1}]
\end{equation*}%
fulfils the system 
\begin{equation}
	\left\{ 
	\begin{array}{ll}
		\begin{array}{l}
			\displaystyle\frac{\partial \upsilon (\widetilde{z})}{\partial t}-\nu
			\triangle \widetilde{z}+\mathrm{curl}\,\upsilon (\widetilde{z})\times y+%
			\mathrm{curl}\,\upsilon (y)\times \widetilde{z}+\nabla \widetilde{\pi } \\ 
			\;\;\qquad \quad \qquad \qquad \qquad =-\mathrm{curl}\,\upsilon
			(z_{1})\times z_{2}-\mathrm{curl}\,\upsilon (z_{2})\times z_{1},%
		\end{array}
		&  \\ 
		\nabla \cdot \widetilde{z}=0 & \mbox{in}\ Q, \\ 
		\widetilde{z}\cdot n=0,\quad (n\cdot D\widetilde{z})\cdot \tau =0 & \mbox{on}%
		\ \Sigma , \\ 
		\widetilde{z}(0)=0 & \mbox{in}\ \mathcal{O}.%
	\end{array}%
	\right.   \label{2dlinearized}
\end{equation}%

\begin{definition}
	The control problem $(\mathcal{P})$ is said to satisfies the second-order
	optimality condition at an optimal pair $(u,S(u))$, if there exists $\delta
	>0,$ such that the coercivity condition 
	\begin{equation}
		J^{\prime \prime }(u)[w,w]>\delta \Vert w\Vert _{L^{2}(Q)}^{2},\quad \forall
		w\in \mathcal{U}_{ad},  \label{coerc}
	\end{equation}%
	\bigskip holds.
\end{definition}

\begin{theorem}
	\label{13:12:20:13:13} Assume the hypothesis of Theorem \ref{main_existence}%
	. Then the control problem $(\mathcal{P})$ satisfies the 
	second-order
	optimality condition at the optimal pair $(u,S(u))$ under the assumption 
	\begin{equation}
		\lambda >2\widehat{K}\lambda _{3}^{2}\lambda _{4},  \label{13:12:20:14:50}
	\end{equation}%
	where $\widehat{K},\lambda _{3},\lambda _{4}$ are defined in Lemmas \ref{L1}%
	, \ref{L4} and \ref{L5}.
\end{theorem}

\textbf{Proof. } Denoting $J(u)=J(u,S(u))$, we have 
\begin{equation*}
	J^{\prime }\left( u\right) [w_{1}]=\int_{0}^{T}\left( \left(
	z_{1},y-y_{d}\right) +\lambda (u,w_{1})\right) dt\ 
\end{equation*}%
by Proposition 4.5 in \cite{AC1}. By similar calculations, made in the
mentioned proposition, we can verify that%
\begin{equation}
	J^{\prime \prime }\left( u\right) [w_{1},w_{2}]=\int_{0}^{T}\left( \left(
	z_{1},z_{2}\right) +\left( \widetilde{z},y-y_{d}\right) +\lambda
	(w_{2},w_{1})\right) dt.  \notag
\end{equation}%
The system \eqref{2dlinearized}, and the duality relation proved in \cite{AC1}, Proposition 5.5, yield   
\begin{equation*}
	J^{\prime \prime }\left( u\right) [w_{1},w_{2}]=\int_{0}^{T}\left( \left(
	z_{1},z_{2}\right) -\left( p,\mathrm{curl}\,\upsilon (z_{1})\times z_{2}+%
	\mathrm{curl}\,\upsilon (z_{2})\times z_{1}\right) +\lambda
	(w_{2},w_{1})\right) dt.
\end{equation*}%
Hence taking $w_{1}=w_{2},$ we have   
\begin{equation*}
	J^{\prime \prime }\left( u\right) [w_{1},w_{1}]=\int_{0}^{T}\left( \Vert
	z_{1}\Vert _{2}^{2}+\lambda \Vert w_{1}\Vert _{2}^{2}-2\left( p,\mathrm{curl}%
	\,\upsilon (z_{1})\times z_{1}\right) \right) dt.
\end{equation*}%
By Lemmas  \ref{L1}, \ref{L4} and \ref{L5}, we deduce the estimate
\begin{align*}
	\int_{0}^{T}\left\vert \left( p,\mathrm{curl}\,\upsilon (z_{1})\times
	z_{1}\right) \right\vert dt& \leq \widehat{K}\int_{0}^{T}\Vert z_{1}\Vert
	_{H^{2}}^{2}\Vert p\Vert _{H^{2}}dt \\
	& \leq \widehat{K}\lambda _{4}\int_{0}^{T}\Vert z_{1}\Vert
	_{H^{2}}^{2}dt\leq \widehat{K}\lambda _{3}^{2}\lambda _{4}\int_{0}^{T}\Vert
	w_{1}\Vert _{2}^{2}dt.
\end{align*}%
Therefore, we conclude that
\begin{equation*}
	J^{\prime \prime }\left( u\right) [w_{1},w_{1}]>(\lambda -2\widehat{K}%
	\lambda _{3}^{2}\lambda _{4})\int_{0}^{T}\Vert w_{1}\Vert _{2}^{2}dt,
\end{equation*}%
which implies the result of the theorem.\hfill $\blacksquare $

\subsection{Global uniqueness of the optimal solution}
Now, we are able to show that the solution of the coupled system 
is unique, and provides the unique global optimal solution for the non-convex optimal control problem  $(\mathcal{P})$.

\begin{theorem}
	\label{21:47:25S} For any%
	\begin{equation*}
		\lambda >2\widehat{K}\lambda _{2}^{2}\lambda _{4},
	\end{equation*}%
	the optimal control problem $(\mathcal{P})$ has a unique global solution,
	where $\widehat{K},$ $\lambda _{3},\lambda _{4}$ are defined in Lemmas \ref%
	{L1}, \ref{L3} and \ref{L5}.
\end{theorem}

\textbf{Proof. } We assume that $u_{1}$ and $u_{2}$ are two optimal control variables
for the problem $(\mathcal{P})$. Let $y_{i}=S(u_{i})$, $i=1,2$, be the
corresponding optimal states with the adjoint states $p_{i}=S^{\ast }(u_{i})$%
, $i=1,2$. Let us set $\overline{u}=u_{2}-u_{1}$, $\overline{\pi }=\pi
_{2}-\pi _{1}$. The differences $\overline{y}=y_{2}-y_{1}$ and $\overline{u}%
=u_{2}-u_{1}$\ verify the system%
\begin{equation}
	\left\{ 
	\begin{array}{ll}
		\begin{array}{l}
			\frac{\partial }{\partial t}\upsilon (\overline{y})=\nu \Delta \overline{y}-%
			\mathrm{curl}\,\upsilon (\overline{y})\times y_{1}-\mathrm{curl}\,\upsilon
			(y_{1})\times \overline{y} \\ 
			\quad \qquad \quad \qquad \qquad \qquad \qquad \qquad -\mathrm{curl}%
			\,\upsilon (\overline{y})\times \overline{y}-\nabla \overline{\pi }+%
			\overline{u}\vspace{1mm},%
		\end{array}
		& \, \\ 
		\mathrm{div}\,\overline{y}=0 & \mbox{in}\ Q, \\ 
		\overline{y}\cdot \mathrm{n}=0,\qquad (\mathrm{n}\cdot D\overline{y})\cdot 
		\mathrm{\tau }=0 & \mbox{on}\ \Sigma , \\ 
		\overline{y}(0)=0 & \mbox{in}\ \mathcal{O}.%
	\end{array}%
	\right.   \label{13:14:28S}
\end{equation}%
The functions $y_{1},\ y_{2}$ are the weak solutions of $(\ref{equation_etat}%
),$ which satisfies $(\ref{var_form_state}).$ Therefore, considering the test function $\phi =p_{1}$ in $(\ref{var_form_state}),$ 
we show that $\overline{y}$ verifies the following variational
equality
\begin{eqnarray*}
	( \tfrac{\partial \overline{y}(t)}{\partial t},p_{1}) _{V}
	&=&-2\nu \left( Dp_{1}(t),D\overline{y}\right) -b\left( p_{1},y_{1},\upsilon
	(\overline{y})\right) +b\left( y_{1},p_{1},\upsilon (\overline{y})\right)  \\
	&&-b\left( p_{1},\overline{y},\upsilon (y_{1})\right) +b\left( \overline{y},p_{1},\upsilon (y_{1})\right) -\left( \mathrm{curl}%
	\,\upsilon (\overline{y})\times \overline{y},p_{1}\right) +(\overline{u}%
	,p_{1}).
\end{eqnarray*}

Writing the adjoint system $(\ref{adj_opt_eq_alpha})$ for $p_{1}$ in a
respective variational form with the test function $\phi =\overline{y}$, we get the equality 
\begin{align*}
	\left( \tfrac{\partial p_{1}}{\partial t},\overline{y}(t)\right) _{V}=& 2\nu
	\left( Dp_{1}(t),D\overline{y}\right) -b\left( \overline{y},p_{1},\upsilon
	(y_{1})\right) +b\left( p_{1},\overline{y},\upsilon (y_{1})\right)  \\
	& +b\left( p_{1},y_{1},\upsilon (\overline{y})\right) -b\left(
	y_{1},p_{1},\upsilon (\overline{y})\right) -\left( y_{1}-y_{d},\overline{y}%
	\right) .
\end{align*}%
Therefore summing the last two equalities, we obtain 
\begin{equation*}
	\frac{\partial }{\partial t}\left( \overline{y}(t),p_{1}\right) _{V}=-(%
	\mathrm{curl}\,\upsilon (\overline{y})\times \overline{y}%
	,p_{1})-(y_{1}-y_{d},\overline{y})+(\overline{u},p_{1}).
\end{equation*}%
The integration over $t\in (0,T)$ and the initial and final conditions for $%
\overline{y}$ and $p_{1}$ give 
\begin{equation}
	0=\int_{0}^{T}-\left( \mathrm{curl}\,\upsilon (\overline{y})\times \overline{%
		y},p_{1}\right) \,-\left( y_{1}-y_{d},\overline{y}\right) +\left( \overline{u%
	},p_{1}\right) \ dt.  \label{1s}
\end{equation}

By the symmetry, we can easily verify that the difference $y_{1}-y_{2}=-
\overline{y}$ satisfies the system%
\begin{equation*}
	\left\{ 
	\begin{array}{ll}
		-\frac{\partial }{\partial t}\upsilon (\overline{y})=-\nu \Delta \overline{y}%
		+\mathrm{curl}\,\upsilon (\overline{y})\times y_{2}+\mathrm{curl}\,\upsilon
		(y_{2})\times \overline{y}-\mathrm{curl}\,\upsilon (\overline{y})\times 
		\overline{y}+\nabla \overline{\pi }-\overline{u}, \\ 
		\mathrm{div}\,\overline{y}=0 \qquad\qquad \qquad\qquad\qquad \qquad \qquad\qquad \qquad \mbox{in}\ Q, \\ 
		\overline{y}\cdot \mathrm{n}=0,\qquad (\mathrm{n}\cdot D\overline{y})\cdot 
		\mathrm{\tau }=0 \qquad\qquad \qquad  \qquad \quad \mbox{on}\ \Sigma , \\ 
		\overline{y}(0)=0\qquad\qquad \qquad\qquad\qquad \qquad \qquad\qquad \qquad \mbox{in}\ \mathcal{O}.%
	\end{array}%
	\right.
\end{equation*}%
Using the same reasoning as above, applied for $\overline{y}$, we deduce
that $y_{1}-y_{2}$ verifies the following relation
\begin{equation}
	0=\int_{0}^{T}-\left( \mathrm{curl}\,\upsilon (\overline{y})\times \overline{%
		y},p_{2}\right) \,+\left( y_{2}-y_{d},\overline{y}\right) -\left( \overline{u%
	},p_{2}\right) \ dt.  \label{2s}
\end{equation}

The sum of the equalities \eqref{1s} and \eqref{2s} yields the identity 
\begin{equation}
	\int_{0}^{T}\Vert \overline{y}\Vert _{2}^{2}\,dt+\int_{0}^{T}\left( 
	\overline{u},p_{1}-p_{2}\right) dt=\int_{0}^{T}\left( \mathrm{curl}%
	\,\upsilon (\overline{y})\times \overline{y},p_{1}+p_{2}\right) \,dt.
	\label{22:17:24S}
\end{equation}%
The optimality condition \eqref{opt1} for the optimal control $u_{1}$ with $%
\Psi =u_{2}$ reads as 
\begin{equation*}
	\int_{0}^{T}\left\{ \lambda (u_{1},\overline{u})+(\overline{u}%
	,p_{1})\,\right\} \,dt\geq 0.
\end{equation*}%
Analogously, the optimality condition \eqref{opt1} 
for the optimal control $%
u_{2}$ with $\Psi =u_{1}$ yields
\begin{equation*}
	-\int_{0}^{T}\left\{ (\lambda u_{2},\overline{u})+(\overline{u}%
	,p_{2})\,\right\} \,dt\geq 0.
\end{equation*}%
Adding the last two inequalities, we deduce 
\begin{equation}
	\lambda \int_{0}^{T}\Vert \overline{u}\Vert _{2}^{2}\,dt\leq \int_{0}^{T}(%
	\overline{u},p_{1}-p_{2})\,dt.  \label{opt120_1}
\end{equation}%
Introducing this relation in \eqref{22:17:24S}, we obtain 
\begin{equation*}
	\lambda \int_{0}^{T}\Vert \overline{u}\Vert _{2}^{2}\,dt+\int_{0}^{T}\Vert 
	\overline{y}\Vert _{2}^{2}\,dt\leq \int_{0}^{T}\left( \mathrm{curl}%
	\,\upsilon (\overline{y})\times \overline{y},p_{1}+p_{2}\right) \,dt.
\end{equation*}%
By Lemmas \ref{L1}, \ref{L3} and \ref{L5} we have 
\begin{align*}
	\int_{0}^{T}\left\vert \left( \mathrm{curl}\,\upsilon (\overline{y})\times 
	\overline{y},p_{1}+p_{2}\right) \right\vert dt& \leq \int_{0}^{T}\Vert 
	\overline{y}\Vert _{H^{2}}^{2}(\Vert p_{1}\Vert _{H^{2}}+\Vert p_{2}\Vert
	_{H^{2}})dt \\
	& \leq 2\widehat{K}\lambda _{4}\int_{0}^{T}\Vert \overline{y}\Vert
	_{H^{2}}^{2}dt\leq 2\widehat{K}\lambda _{2}^{2}\lambda _{4}\int_{0}^{T}\Vert 
	\overline{u}\Vert _{2}^{2}dt.
\end{align*}%
Therefore 
\begin{equation*}
	(\lambda -2\widehat{K}\lambda _{2}^{2}\lambda _{4})\int_{0}^{T}\Vert 
	\overline{u}\Vert _{2}^{2}+\int_{0}^{T}\Vert \overline{y}\Vert _{2}^{2}\
	dt\leq 0,
\end{equation*}%
which implies the result of the theorem.\hfill $\blacksquare $

\bigskip

\bigskip

	Let us remark that in  Theorem \ref{21:47:25S}, the constants $\lambda _{2},$ $%
	\lambda _{4}$  just  depend 
	 on the bounded set 
	 $\mathcal{U}_{ad}$ of the 
	 admissible control variables and the target and initial velocities $y_{d},$ $y_{0}$.
	Hence, from pratical point of view,  the  motion of 
	a second grade fluid
	can be optimally controlled by taking a sufficiently large intensity of the cost $\lambda.$

\bigskip

\subsection*{Acknowledgments}
The research of Nikolai V. Chemetov was supported by the project n. 2021 / 03758-8 of Regular Research Grant - FAPESP, Brazil and by the Program of "Docentes Novos", Process no: 20.1.4175.1.0. of Universidade de São Paulo.

\medskip

The work of F. Cipriano was partially supported by the Funda\c{c}\~{a}o para a Ci\^{e}ncia e a Tecnologia (Portuguese Foundation for Science and
Technology) through the project UID/MAT/00297/2019 (Centro de Matem\'{a}tica e Aplica\c{c}\~{o}es).

\end{document}